\documentclass[a4paper,10pt]{amsart}
\usepackage{amsmath}
\usepackage{amsfonts}
\usepackage{amssymb}
\usepackage[cp1250]{inputenc}

\newtheorem{theorem}{Theorem}[section]
\newtheorem{proposition}[theorem]{Proposition}

\newtheorem{definition}[theorem]{Definition}

\newcommand{\dow}[1][]{\textbf{Proof{#1}. }}
\newcommand{\mdeg}{\mbox{mdeg}\,}
\newcommand{\wdeg}{\mbox{deg}_w\,}
\newcommand{\wmdeg}{\mbox{mdeg}_w\,}
\newcommand{\Aut}{\mbox{Aut}\,}
\newcommand{\Tame}{\mbox{Tame}\,}
\newcommand{\length}{\mbox{length}\,}

\begin{document}
\author{Marek Kara\'{s}}
\title{On weighted bidegree of polynomial automorphisms of $\mathbb{C}^2$}
\date{}
\maketitle

%\footnote{keywords{polynomial automorphism, tame automorphism, wild automorphism, multidegree.\\
%\textit{2010 Mathematics Subject Classification:} 14Rxx,14R10}}

\begin{abstract}
Let $F=(F_1,F_2):\mathbb{C}^2\rightarrow\mathbb{C}^2$ be a polynomial automorphism. It is well know that $\deg F_1 | \deg F_2$ or $\deg F_2 | \deg F_1.$ On the other hand, if $(d_1,d_2)\in\mathbb{N}_+^2=(\mathbb{N}\setminus\{ 0 \} )^2$ is such that $d_1 | d_2$ or $d_2 | d_1,$ then one can construct a polynomial automorphism $F=(F_1,F_2)$ of $\mathbb{C}^2$ with $\deg F_1=d_1$ and $\deg F_2=d_2.$ 

Let us fix $w=(w_1,w_2)\in\mathbb{N}_+^2$ and consider the \textit{weighted degree} on $\mathbb{C}[x,y]$ with $\wdeg x=w_1$ and $\wdeg y=w_2.$ In this note we address the structure of the set $\{ (\wdeg F_1,\wdeg F_2) \, : \allowbreak\, (F_1,F_2) \mbox{ is an}\allowbreak \mbox{automorphism of } \mathbb{C}^2\,\} .$
\end{abstract}

\section{Introduction}

Let us fix $n$-tuple $w=(w_1,\ldots,w_n)\in\mathbb{N}_+^n=(\mathbb{N}\setminus \{ 0\})^n.$ In this note we will write $\deg h$ for the usal total degreee of a polynomia $h\in\mathbb{C}[x_1,\ldots,x_n]$ and $\wdeg h$ for the \textit{weighted degree} of $h$ with respect to $w,$ where  
\begin{equation}
\deg_w h=\max \left\{ \alpha_1 w_1+\cdots +\alpha_n w_n\,: c_{\alpha}\neq 0\, \right\}
\end{equation}  
for
\begin{equation}
h=\sum_{\alpha=(\alpha_1,\ldots,\alpha_n)\in\mathbb{N}^n}c_{\alpha}x_1^{\alpha_1}\cdots x_n^{\alpha_n}.%\in\mathbb{C}[x_1,\ldots,x_n].
\end{equation}
In other words, we assume that $\wdeg x_1=w_1,\ldots,\wdeg x_n=w_n.$

If $F=(f_{1},\ldots,f_{n}):\mathbb{C}^n\rightarrow\mathbb{C}^n$ is a polynomial mapping, then by \textit{multidegree} of $F$ we mean the following $n$-tuple $\mdeg  F=(\deg f_{1},\ldots,\allowbreak\deg f_{n})\in\mathbb{N}^n$ and by the \textit{weighted multidegree} of $F$ with respect to the \textit{weight} $w$ we mean the following one  $\wmdeg  F=(\deg_w f_{1},\ldots,\allowbreak\deg_w f_{n})\in\mathbb{N}^n.$  Sometimes, when $n=2,$ $\mdeg F$ is also called \textit{bidegree} of $F$ and $\wmdeg F$ is called \textit{weighted bidegree} of $F.$ 

Let us recall that a polynomial automorphism $F$ of $\mathbb{C}^n$ is called \textit{tame} if it can be obtained as a composition of affine and triangular automorphisms. As usual a mapping $G=(G_1,\ldots,G_n):\mathbb{C}^n\rightarrow\mathbb{C}^n$ is called \textit{affine} if $\deg G_i=1$ for $i=1,\ldots n,$ and a mapping $H=(H_1,\ldots,H_n):\mathbb{C}^n\rightarrow\mathbb{C}^n$ is called \textit{triangular} if for some permutation $\sigma$ of $\{ 1,\ldots,n\}$ we have $H_{\sigma(i)}=c_i\cdot x_{\sigma(i)}+h_i$  for $i=1,\ldots,n$ and some $c_i\in\mathbb{C}^*=\mathbb{C}\setminus\{ 0\},$ $h_i\in\mathbb{C}[x_{\sigma(1)},\ldots,x_{\sigma(i-1)}].$

In what follows we will write $\Aut (\mathbb{C}^{n})$ for the group of the all
polynomial automorphisms of $\mathbb{C}^{n}$ and $\Tame (\mathbb{C}^n)$ for the subgroup of $\Aut (\mathbb{C}^{n})$ containing all the tame automorphisms. 
 Then, one can consider two functions
(also denoted $\mdeg$ and $\wmdeg$) mapping $\Aut (\mathbb{C}^{n})$ into $%
\mathbb{N}_{+}^{n}.$ 
It is well-known \cite{Jung,Kulk} that 
\begin{equation}\label{row_jungKulk}
\mdeg(\Aut(\mathbb{C}^2))=\mdeg(\Tame(\mathbb{C}^2))=\{ (d_1,d_2)\in\mathbb{N}_+^2 \,:\, d_1 | d_2\mbox{ or } d_2 | d_1 \}.
\end{equation}
Since $\Aut(\mathbb{C}^2)=\Tame(\mathbb{C}^2),$ we obviously have
\begin{equation}
\wmdeg(\Aut(\mathbb{C}^2))=\wmdeg(\Tame(\mathbb{C}^2)).
\end{equation}

This note address the structure of the above set. Namely we show the following

\begin{theorem}\label{Tw_main}
Let $w=(w_1,w_2)\in\mathbb{N}_+^2.$ Then the set $\wmdeg (\Aut (\mathbb{C}^2))$ is equal to
\begin{eqnarray*}
&& \left\{ (d_1,d_2)\in(w_1\mathbb{N}_+)^2 \, : \, d_1|d_2\mbox{ or } d_2|d_1 ,\, \max \{ d_1,d_2\} \geq \widetilde{w}, \, 
\min \{ d_1,d_2 \}< \widetilde{w} \Rightarrow %\mbox{ or }
\min\{ d_1,d_2 \}=\underline{w} \, \right\} \\
&\cup&  \left\{ (d_1,d_2)\in(w_2\mathbb{N}_+)^2 \, :\,  d_1|d_2\mbox{ or } d_2|d_1 ,\, \max \{ d_1,d_2\} \geq \widetilde{w}, \, 
\min \{ d_1,d_2 \}< \widetilde{w} \Rightarrow %\mbox{ or }
\min\{ d_1,d_2 \}=\underline{w} \, \right\} \\
&\cup &  \left\{ (w_1,w_2), (w_2,w_1), (\tilde{w},\tilde{w} )\right\} ,
\end{eqnarray*}
where $\widetilde{w}:=\max \{ w_1,w_2 \}$ and $\underline{w}:=\min\{ w_1,w_2\} .$
\end{theorem}

Notice that if $w_1=w_2,$ then the set given on the right-hand side of the above equality is equal to $\{ (d_1,d_2)\in(w_1\mathbb{N}_+)^2 \, :\, d_1|d_2\mbox{ or } d_2|d_1\, \} .$ In particular, for $(w_1,w_2)=(1,1),$ one obtain the equality (\ref{row_jungKulk}).

For information about multidegrees of tame and wild automorphisms of $\mathbb{C}^3$ see \cite{Jiantao&Xiankun, Karas1,Karas2,Karas3,Karas4,Karas5,KarasZygad, KarasZygad2,Zygadlo}. 
%--------------------------------------------------------------------------------
\section{Lenght of $F\in\Aut(\mathbb{C}^2)$ and the weighted bidegree}

In this section we show that $\wmdeg(\Aut(\mathbb{C}^2))$ is contained in the set given on the right-hand side of the equality of Theorem \ref{Tw_main}. More precisely we show Theorems \ref{tw_lenght_1} and \ref{Tw_lenght_2} below, but we start with the following 

\begin{proposition}[see e.g. {\cite[Prop. 9.2]{Karas4}}]
\label{Prop_trojk_minim}
If $F\in \Aut\left( \mathbb{C}^{2}\right) ,$
then there is a number $l\in \mathbb{N}$ (possibly zero), affine automorphisms $%
L_{1},L_{2}$ of $\mathbb{C}^{2}$ and triangular automorphisms $T_{1},\ldots
,T_{l}$ of the forms 
\begin{align}
T_{i} &:\mathbb{C}^{2}\ni \left( x,y\right) \mapsto \left( x,y+f_{i}(x)\right)
\in \mathbb{C}^{2}\qquad \text{for }i=1,3,\ldots ,  \label{Row_Prop_amalg_1} \\
T_{i} &:\mathbb{C}^{2}\ni \left( x,y\right) \mapsto \left( x+f_{i}(y),y\right)
\in \mathbb{C}^{2}\qquad \text{for }i=2,4,\ldots ,  \label{Row_Prop_amalg_2}
\end{align}
with $\deg f_{i}>1,$ such that 
\begin{equation*}
F=L_{2}\circ T_{l}\circ \cdots \circ T_{1}\circ L_{1}.
\end{equation*}
Moreover, the number $l$ is unique, and one can require that $T_{i},$ $%
i=1,\ldots ,l,$ are of the form (\ref{Row_Prop_amalg_1}) for even $i$ and of
the form (\ref{Row_Prop_amalg_2}) for odd $i.$
\end{proposition}

\begin{definition}[see e.g. {\cite[p.612]{Furter}}]
Let $F\in \Aut\left( \mathbb{C%
}^{2}\right) $ be a polynomial automorphism. The number $l$ from Proposition 
\ref{Prop_trojk_minim} is called the length of $F$ and denoted $\length F.$
\end{definition}

Now, we are in a position to prove Theorems \ref{tw_lenght_1} and \ref{Tw_lenght_2}.

\begin{theorem}\label{tw_lenght_1}
Let us fix $w=(w_1,w_2)\in\mathbb{N}_+^2.$
If $F:\mathbb{C}^2\rightarrow\mathbb{C}^2$ is a polynomial automrphism with $\length F\leq 1,$ then the weighted multidegree $\wmdeg F$ is an element of the following set
\begin{eqnarray*}
&& \left\{ (w_1,kw_1),\, (kw_1,w_1),\, (kw_1,kw_1) : k\in\mathbb{N}_+\mbox{ and } kw_1\geq w_2 \right\} \\
&\cup&  \left\{ (w_2,kw_2),\, (kw_2,w_2),\, (kw_2, kw_2) : k\in\mathbb{N}_+\mbox{ and } kw_2\geq w_1 \right\}  \\
&\cup&   \left\{ (w_1,w_2), (w_2,w_1), (\tilde{w},\tilde{w} )\right\} ,
\end{eqnarray*}
where $\tilde{w}:=\max \{ w_1,w_2 \}.$
\end{theorem}

\dow If $\length F=0,$ then $F$ is affine and so one can easy check that $\wmdeg F$ belongs to $\left\{ (w_1,w_2), (w_2,w_1),\allowbreak (\tilde{w},\tilde{w} ) \right\} .$

Assume that $\length F=1.$ By Proposition \ref{Prop_trojk_minim} we can assume that $F=L_2\circ T \circ L_1,$ where $L_1,L_2$ are affine automorphisms and $T$ is of the form
\begin{equation}
T:\mathbb{C}^2\ni (x,y)\mapsto (x,y+f(x))\in\mathbb{C}^2,
\end{equation}
with $\deg f>1.$

We have three cases: (I) $\wmdeg L_1=(w_1,w_2),$ (II) $\wmdeg L_1=(w_2,w_1)$ and (III) $\wmdeg L_1=(\tilde{w},\tilde{w}).$ Thus we have
\begin{equation}\label{Row_k1k2}
(k_1,k_2):=\wmdeg (T\circ L_1)=
\left\{
\begin{array}{ll}
(w_1,\max\{ w_1\cdot\deg f\, ,\, w_2 \} ), & \mbox{for case (I),}\\
(w_2,\max\{ w_2\cdot\deg f\, ,\, w_1 \} ), & \mbox{for case (II),}\\
(\tilde{w},\tilde{w}\cdot\deg f), & \mbox{for case (III).}
\end{array}
\right. 
\end{equation}
Since $L_2$ is affine, it follows that $(d_1,d_2):=\wmdeg F=\wmdeg (L_1\circ T\circ L_1)$ belongs to $\{ (k_1,k_2), (k_2,k_1),(\tilde{k},\tilde{k}) \},$ where $\tilde{k}:=\max\{ k_1, k_2 \} .$ Thus, we have:
\begin{list}{}{}

\item \textsc{Case (I).} If $\max\{ w_1\cdot\deg f\, ,\, w_2 \} =w_2,$ then $(k_1,k_2)=(w_1,w_2)$ and so $(d_1,d_2)$ belongs to $\{ (w_1,w_2),\, (w_2,w_1),\, (w_2,w_2) \}$ else $(k_1,k_2)=(w_1,w_1\cdot\deg f)$ and so $(d_1,d_2)$ belongs to $\{ (w_1,w_1\cdot\deg f),\, (w_1\cdot\deg f,w_1),\, (w_1\cdot\deg f,w_1\cdot\deg f) \} .$\vspace{0.2cm}

\item \textsc{Case (II).} If $\max\{ w_2\cdot\deg f\, ,\, w_2 \} =w_1,$ then $(k_1,k_2)=(w_2,w_1)$ and so $(d_1,d_2)$ belongs to $\{ (w_1,w_2),\, (w_2,w_1),\, (w_1,w_1) \}$ else $(k_1,k_2)=(w_2,w_2\cdot\deg f)$ and so $(d_1,d_2)$ belongs to $\{ (w_2,w_2\cdot\deg f),\, (w_2\cdot\deg f,w_2),\, (w_2\cdot\deg f,w_2\cdot\deg f) \} .$\vspace{0.2cm}

\item \textsc{Case (III).} If $\tilde{w}=w_1,$ then $(k_1,k_2)=(w_1,w_1\cdot\deg f)$ and so  $(d_1,d_2)$ belongs to $\{ (w_1,w_1\cdot\deg f),\, (w_1\cdot\deg f,w_1),\, (w_1\cdot\deg f,w_1\cdot\deg f) \} $ else $\tilde{w}=w_2,$ $(k_1,k_2)=(w_2,w_2\cdot\deg f)$ and so $(d_1,d_2)$ belongs to $\{ (w_2,w_2\cdot\deg f),\, (w_2\cdot\deg f,w_2),\, (w_2\cdot\deg f,w_2\cdot\deg f) \} .$
\end{list}
Thus the result follows. 
$\Box$

\begin{theorem}\label{Tw_lenght_2}
Let us fix $w=(w_1,w_2)\in\mathbb{N}_+^2.$
If $F:\mathbb{C}^2\rightarrow\mathbb{C}^2$ is a polynomial automrphism with $\length F\geq 2,$ then the weighted multidegree $\wmdeg F$ is an element of the following set
\begin{eqnarray*}
&& \left\{ (d_1,d_2)\in(w_1\mathbb{N}_+)^2 \, : \min \{ d_1,d_2\} \geq \max \{ w_1,w_2\} \mbox{ and } (\, d_1|d_2\mbox{ or } d_2|d_1 \, )\right\} \\
&\cup&  \left\{ (d_1,d_2)\in(w_2\mathbb{N}_+)^2 \, : \min \{ d_1,d_2\} \geq  \max \{ w_1,w_2\} \mbox{ and } (\, d_1|d_2\mbox{ or } d_2|d_1 \, ) \right\} .
\end{eqnarray*}
In particular $| \wmdeg F|>|w|$ when $\length F\geq 2.$ 
\end{theorem}

\dow Let $l:=\length F.$ By Proposition \ref{Prop_trojk_minim} we can assume that 
\begin{equation*}
F=L_{2}\circ T_{l}\circ \cdots \circ T_{1}\circ L_{1},
\end{equation*}
where $L_1,L_2$ are affine automorphism of $\mathbb{C}^2$ and $T_{1},\ldots
,T_{l}$ are triangular automorphisms of the forms 
\begin{align}
T_{i} &:\mathbb{C}^{2}\ni \left( x,y\right) \mapsto \left( x,y+f_{i}(x)\right)
\in \mathbb{C}^{2}\qquad \text{for }i=1,3,\ldots ,   \\
T_{i} &:\mathbb{C}^{2}\ni \left( x,y\right) \mapsto \left( x+f_{i}(y),y\right)
\in \mathbb{C}^{2}\qquad \text{for }i=2,4,\ldots .
\end{align}

Let $(k_1,k_2):=\wmdeg (T_1\circ L_1).$ By (\ref{Row_k1k2}), we have that $k_2>k_1,$ $k_2\geq\max \{ w_1,w_2 \}$ and $k_2\in w_1\mathbb{N}_+\cup w_2\mathbb{N}_+ .$ It is easy to see that
\begin{equation}
\wmdeg (T_2\circ T_1\circ L_1)=(k_2\cdot\deg f_2\, ,\, k_2)
\end{equation}
and for $l>2$
\begin{equation}
\wmdeg (T_l\circ\cdots\circ T_1\circ L_1)=
\left\{
\begin{array}{ll}
(k_2\prod_{i=2}^l\deg f_i\, ,\, k_2\prod_{i=2}^{l-1}\deg f_i), & \mbox{for even }l, \vspace{0.3cm} \\
(k_2\prod_{i=2}^{l-1}\deg f_i\, ,\, k_2\prod_{i=2}^{l}\deg f_i), & \mbox{for odd }l.

\end{array}
\right.
\end{equation}
Thus
\begin{equation}
\wmdeg (L_2\circ T_l\circ\cdots\circ T_1\circ L_1)\in\{ (m_1,m_2),(m_2,m_1),(m_2,m_2) \},
\end{equation}
where $m_1:= k_2\prod_{i=2}^{l-1}\deg f_i$ and $m_2:=k_2\prod_{i=2}^{l}\deg f_i$ for $l>2,$ and $m_1:=k_2$ and $m_2:=k_2\deg f_2$ for $l=2.$ 
Hence, the result follows. $\Box$

%-----------------------------------------------------------------------------
\section{Examples}

Let $Z$ denotes the set given in Theorem \ref{Tw_main}.
By Theorems \ref{tw_lenght_1} and \ref{Tw_lenght_2}, in order to proof Thoerem \ref{Tw_main}, it is enough to show an example of automorphism $F\in\Aut(\mathbb{C}^2)$ with $\wmdeg F=(d_1,d_2)$ for each $(d_1,d_2)\in Z.$  

Without lose of generality we can assume that $w_1\leq w_2.$ First consider the case $w_1=w_2.$ Take any $(d_1,d_2)\in Z.$ Since $w_2|d_1$ and $w_2|d_2,$ it follows that one can take
\begin{equation}
F=\left\{
\begin{array}{ll}
T_2\circ T_1, & \mbox{for } d_1<d_2, \\
\tilde{T}_2\circ\tilde{T}_1, & \mbox{for } d_1>d_2, \\
L\circ T_1, & \mbox{for } d_1=d_2,
\end{array}
\right. 
\end{equation} 
where $T_1(x,y)=(x+y^{\frac{d_1}{w_2}},y),$ $T_2(x,y)=(x,y+x^{\frac{d_2}{d_1}},$ $\tilde{T}_1(x,y)=(x,y+x^{\frac{d_2}{w_2}}),$ $\tilde{T}_2(x,y)=(x+y^{\frac{d_2}{d_1}},y)$ and $L(x,y)=(x,y+x).$

Now, consider the case $w_1<w_2$ and take any $(d_1,d_2)\in Z.$ If $(d_1,d_2)\in(w_2\mathbb{N}_+)^2,$ then one can take
\begin{equation}
F=\left\{
\begin{array}{ll}
T_2\circ T_1, & \mbox{for } d_1<d_2, \\
\tilde{T}_2\circ\tilde{T}_1\circ\tilde{L}, & \mbox{for } d_1>d_2, \\
L\circ T_1, & \mbox{for } d_1=d_2,
\end{array}
\right. 
\end{equation} 
where $T_1,T_2,\tilde{T}_1,\tilde{T}_2$ and $L$ are defined as in the case $w_1=w_2,$ and $\tilde{L}(x,y)=(y,x).$

If $(d_1,d_2)\in(w_1\mathbb{N}_+)^2,$ then we have two cases: (I) $\min \{ d_1,d_2 \}\geq w_2$ and (II) $\min \{ d_1,d_2 \}=w_1.$ In case (I) one can take
\begin{equation}
F=\left\{
\begin{array}{ll}
T_2\circ T_1, & \mbox{for } d_1<d_2, \\
\tilde{T}_2\circ\tilde{T}_1, & \mbox{for } d_1>d_2, \\
L\circ T_1, & \mbox{for } d_1=d_2,
\end{array}
\right. 
\end{equation}
where $T_1(x,y)=(x,y+x^{\frac{d_1}{w_1}}),$ $T_2(x,y)=(x+y^{\frac{d_2}{d_1}}),$ $\tilde{T}_1(x,y)=(x,y+x^{\frac{d_2}{w_1}}),$ $\tilde{T}_2(x,y)=(x+y^{\frac{d_1}{d_2}})$ and $L(x,y)=(x+y,y).$

And, in case (II), we one can take
\begin{equation}
F=\left\{
\begin{array}{ll}
T_1, & \mbox{for } d_1<d_2, \\
\tilde{L}\circ T_2, & \mbox{for } d_1>d_2, \\
L\circ T_1, & \mbox{for } d_1=d_2,
\end{array}
\right. 
\end{equation}
where $T_1(x,y)=(x,y+x^{\frac{d_2}{w_1}}),$ $T_2(x,y)=(x,y+x^{\frac{d_1}{w_1}}),$ $L(x,y)=(x+y,y)$ and  $\tilde{L}(x,y)=(y,x).$

Finally, if $(d_1,d_2)\in\{ (w_1,w_2),(w_2,w_1),(w_2,w_2) \},$ then one can take
\begin{equation}
f(x,y)=\left\{
\begin{array}{ll}
(x,y), & \mbox{for } (d_1,d_2)=(w_1,w_2), \\
(y,x), & \mbox{for } (d_1,d_2)=(w_2,w_1), \\
(x+y,x) & \mbox{for } (d_1,d_2)=(w_2,w_2). \\
\end{array}
\right.
\end{equation}

%-----------------------------------------------------------------------------

\vspace{1cm}

\textsc{Marek Kara\'{s}\newline
Instytut Matematyki,\newline
Wydział Matematyki i Informatyki\newline
Uniwersytetu Jagiello\'{n}skiego\newline
ul. \L ojasiewicza 6}\newline
\textsc{30-348 Krak\'{o}w\newline
Poland\newline
} e-mail: Marek.Karas@im.uj.edu.pl\vspace{0.5cm}

%and\vspace{0.5cm}

%\textsc{Jakub Zygad\l o\newline
%Instytut Informatyki,\newline
%Wydział Matematyki i Informatyki\newline
%Uniwersytetu Jagiello\'{n}skiego\newline
%ul. \L ojasiewicza 6}\newline
%\textsc{30-348 Krak\'{o}w\newline
%Poland\newline
%} e-mail: Jakub.Zygadlo@ii.uj.edu.p

\end{document}